\documentclass[11pt]{article}
\usepackage{amsmath,amssymb}

\def\N{\bf \mbox{I\hspace{-.15em}N}}

\def\R{{\bf \mbox{I\hspace{-.20em}R}}}

\def\C{C^{\infty}(M, {\bf \mbox{I\hspace{-.20em}R}})}

\newtheorem{definition}{Definition}[section]

\newtheorem{proposition}[definition]{Proposition}
\newtheorem{theorem}[definition]{Theorem}

\newtheorem{examples}[definition]{Examples}

\newenvironment{proof}{\noindent{\bf Proof.}}{\hfill $\square$}
\newenvironment{proof-of}{\noindent{\bf Proof of}}{\hfill $\square$}

\begin{document}

\title{Characteristic foliation of twisted Jacobi manifolds}
\author{J. M. Nunes da Costa \\Departamento de Matem\'atica \\ Universidade de Coimbra\\Apartado 3008
\\3001-454 Coimbra - Portugal\\ {\small E-mail: jmcosta@mat.uc.pt}
\and F. Petalidou\\Faculty of Sciences and Technology
\\University of Peloponnese \\22100 Tripoli - Greece\\{\small E-mail: petalido@uop.gr}}

\date{}
\maketitle

\begin{abstract}
We study the characteristic foliation of a twisted Jacobi
manifold. We show that a twisted Jacobi manifold is foliated into
leaves that are, according to the parity of the dimension, endowed
with a twisted contact or a twisted locally conformal symplectic
structure.
\end{abstract}

\vspace{3mm} \noindent {\bf{Keywords : }}{Twisted Jacobi manifold,
characteristic foliation.}

\vspace{3mm} \noindent {\bf A.M.S. classification (2000):}{53C12,
53D10, 53D17.}

\section{Introduction}

The notion of twisted Jacobi manifold was introduced by the
authors in \cite{jf} and several of its properties and relations
with other geometric structures, such as twisted Dirac-Jacobi
bundles and quasi-Jacobi bialgebroids, were studied in
\cite{jf-tj}. Twisted Jacobi manifolds appeared as a natural
generalization of the twisted Poisson manifolds, which were
introduced by Severa and Weinstein in \cite{sw}, motivated by
works on topological field theory \cite{p} and on string theory
\cite{kl}. In \cite{jf-tj}, some examples of twisted Jacobi
structures on manifolds were presented, including twisted locally
conformal symplectic structures. In this Note, we show that
twisted contact structures
 also provide examples of twisted Jacobi structures. Twisted contact and twisted locally conformal symplectic structures are
 two important types of twisted Jacobi structures on a manifold. In
 fact, we prove that, according to the parity of its dimension, a
 transitive twisted Jacobi manifold is either a twisted contact
 manifold, or a twisted locally conformal symplectic manifold.
 The characteristic foliation of a twisted Jacobi manifold is also
 discussed in this Note and we show that each characteristic leaf of
a twisted Jacobi manifold is endowed with a transitive twisted
Jacobi structure.

The paper starts with a very brief review, in section 2, of the
main properties of twisted Jacobi manifolds. Section 3 is devoted
to the study of the characteristic foliation of a twisted Jacobi
manifold.


\section{Twisted Jacobi manifolds}
A \emph{twisted Jacobi manifold} (\cite{jf, jf-tj}) is a
differentiable manifold $M$ equipped with a bivector field
$\Lambda$, a vector field $E$ and a $2$-form $\omega$ such that
\begin{equation}\label{def-tj}
\frac{1}{2}[(\Lambda,E),(\Lambda,E)]^{(0,1)}=(\Lambda,E)^{\#}(d\omega,\omega).
\end{equation}
In (\ref{def-tj}), $[\cdot,\cdot]^{(0,1)}$ denotes the Schouten
bracket of the Lie algebroid $(TM\times \R, [\cdot,\cdot], \pi)$
over $M$, modified by the 1-cocycle $(0,1)$ of the Lie algebroid
cohomology complex with trivial coefficients \cite{im}, and
$(\Lambda,E)^{\#}$ is the natural extension of $(\Lambda,E)^{\#} :
\Gamma(T^*M\times \R)\to \Gamma(TM\times \R)$, given, for all
$(\alpha,f)\in \Gamma(T^*M\times \R)$, by
$(\Lambda,E)^{\#}(\alpha,f)=(\Lambda^{\#}(\alpha)+fE,-\langle
\alpha,E\rangle)$, to a homomorphism from
$\Gamma(\bigwedge^k(T^*M\times \R))$ to
$\Gamma(\bigwedge^k(TM\times \R))$, $k\in \N$, ($k=3$, in
(\ref{def-tj})) defined, for all $(\eta,\xi) \in \Gamma
(\bigwedge^k( T^*M \times \R))$ and $(\alpha_1,f_1),
\ldots,(\alpha_k,f_k)\in \Gamma(T^*M \times \R)$, by
\begin{eqnarray*}
\lefteqn{(\Lambda,E)^\#(\eta,\xi)((\alpha_1,f_1), \ldots,
(\alpha_k,f_k)) }  \nonumber \\
& & =(-1)^k(\eta,\xi)((\Lambda,E)^\# (\alpha_1,f_1), \cdots ,
(\Lambda,E)^\# (\alpha_k,f_k)).
\end{eqnarray*}

For a bivector field $\Lambda$ on $M$, we consider the usual
homomorphism   $\Lambda^{\#}: \Gamma(T^*M) \to \Gamma(TM)$
associated to $\Lambda$ and we define its natural extension
$\Lambda^{\#}: \Gamma(\bigwedge^k T^*M) \to \Gamma(\bigwedge^k
TM)$, $k \in \N$, by setting, for all $\eta \in \Gamma
(\bigwedge^k T^*M )$ and $\alpha_1, \ldots, \alpha_k \in
\Gamma(T^*M)$,
$$
\Lambda^\#(\eta)(\alpha_1, \ldots, \alpha_k )= (-1)^k \eta
(\Lambda^\#(\alpha_1), \ldots, \Lambda^\#(\alpha_k)).
$$
Also, following \cite{sw}, we denote by $(\Lambda^{\#}\otimes
1)(\eta)$ the section of $(\bigwedge^{k-1}TM)\otimes T^*M$ that
acts on multivector fields by contraction with the factor in
$T^*M$: for all $X\in \Gamma(TM)$ and $\alpha_1,\ldots,
\alpha_{k-1} \in \Gamma(T^*M)$,
$$
(\Lambda^{\#}\otimes 1)(\eta)(\alpha_1,\ldots, \alpha_{k-1})(X) =
(-1)^k\eta(\Lambda^{\#}(\alpha_1),\ldots,\Lambda^{\#}(\alpha_{k-1}),X).
$$

The next Proposition gives an equivalent expression of
(\ref{def-tj}) in terms of the usual Schouten bracket.

\begin{proposition}(\cite{jf-tj})\label{eq-def}
The pair $((\Lambda,E),\omega)$, with $(\Lambda,E) \in
\Gamma(\bigwedge^2(TM \times \R))$ and $\omega \in
\Gamma(\bigwedge^2T^*M)$, defines a  twisted Jacobi structure on
$M$ if and only if
\begin{equation}\label{eq-def-1}
[\Lambda,\Lambda] +2E \wedge \Lambda  =  2 \Lambda^{\#}(d\omega)+
2( \Lambda^{\#} \omega)\wedge E
\end{equation}
and
\begin{equation} \label{eq-def-2}
[E,\Lambda]=  ( \Lambda^{\#} \otimes 1)(d\omega)(E)- ((
\Lambda^{\#} \otimes 1)(\omega)(E)) \wedge E.
\end{equation}
\end{proposition}

\vspace{3mm}

As in the case of a Jacobi manifold \cite{lch}, given a twisted
Jacobi structure $((\Lambda,E),\omega)$ on $M$, $(\Lambda,E)$
defines on $\C$ the internal composition law
\begin{equation}\label{croch-fun}
\{f,g\}=\Lambda(df,dg)+\langle fdg-gdf,E\rangle,
\hspace{5mm}f,g\in \C,
\end{equation}
that is bilinear, skew-symmetric but it does not, in general,
satisfy the Jacobi identity. We have \cite{jf-tj}, for all $f,g,h
\in  \C$,
$$
\{f,\{g,h\}\} + \{g,\{h,f\}\} + \{h,\{f,g\}\}
=(\Lambda,E)^{\#}(d\omega,\omega)((df,f),(dg,g),(dh,h)).
$$
Therefore, $(\C,\{\cdot,\cdot\})$ is not a local Lie algebra.
However, $((\Lambda,E),\omega)$ defines a Lie algebroid structure
$(\{\cdot,\cdot\}^{\omega}, \pi \circ (\Lambda,E)^\#)$ on the
vector bundle $T^*M\times \R \to M$, \cite{jf-tj}. The bracket on
the space $\Gamma(T^*M\times \R)$ of smooth sections is given, for
all $(\alpha,f), (\beta,g)\in \Gamma(T^*M\times \R)$, by
\begin{equation*}\label{croch-form}
\{ (\alpha,f),(\beta,g) \} ^{\omega} = \{(\alpha,f),(\beta,g)\} +
(d \omega,\omega)((\Lambda,E)^\# (\alpha,f), (\Lambda,E)^\#
(\beta,g), \cdot),
\end{equation*}
where $\{\cdot,\cdot\}$ denotes the Kerbrat-Souici-Benhammadi
bracket (\cite{krb}) and the anchor map is $\pi \circ
(\Lambda,E)^\#$, where $\pi :TM\times \R \to TM$ denotes the
projection on the first factor.

\vspace{2mm}

Next, we present two important examples of twisted Jacobi
structures on a manifold.
\begin{examples}
\end{examples}
\vspace{-2mm}

\noindent \emph{1. Twisted locally conformal symplectic
manifolds:} A \emph{twisted locally conformal symplectic manifold}
\cite{jf-tj} is a manifold $M$ of even dimension $2n$ equipped
with a non-degenerate $2$-form $\Theta$, a closed $1$-form
$\vartheta$,
and a $2$-form $\omega$ such that
$$
d(\Theta - \omega) + \vartheta \wedge (\Theta - \omega) =0.
$$
Let $E$ be the unique vector field and $\Lambda$ the unique
bivector field on $M$ which are defined by
\begin{equation}\label{eq-tlcs}
i(E)\Theta = -\vartheta  \hspace{5mm} \mathrm{and} \hspace{5mm}
i(\Lambda^\#(\alpha))\Theta = -\alpha, \hspace{5mm}
\mathrm{for}\,\, \mathrm{all}\,\, \alpha \in \Gamma(T^*M).
\end{equation}
Then, we have
$$
E = \Lambda^\#(\vartheta) \hspace{3mm} \mathrm{and} \hspace{3mm}
\Lambda = \Lambda^\#(\Theta).
$$
By a simple, but very long computation, we prove that the pair
$((\Lambda,E),\omega)$ satisfies the relations (\ref{eq-def-1})
and (\ref{eq-def-2}). Whence, $((\Lambda,E),\omega)$ endows $M$
with a twisted Jacobi structure.

\vspace{1mm}

\noindent \emph{2. Twisted contact manifolds:} A \emph{twisted
contact manifold} is a manifold $M$ of odd dimension $2n+1$
equipped with a $1$-form $\vartheta$ and a $2$-form $\omega$ such
that $\vartheta\wedge (d\vartheta + \omega)^n\neq 0$, everywhere
in $M$. Let us consider on $M$ the vector field $E$ defined by
$$
i(E)\vartheta =1 \;\;\; \mathrm{and} \;\;\; i(E)(d\vartheta
+\omega)=0,
$$
and the bivector field $\Lambda$ whose associated morphism
$\Lambda^{\#}$ is given, for all $\alpha\in\Gamma(T^*M)$, by
$$
\Lambda^{\#}(\vartheta)=0\;\;\; \mathrm{and}
\;\;\;i(\Lambda^{\#}(\alpha))(d\vartheta + \omega) = -(\alpha-
\langle \alpha,E\rangle \, \vartheta).
$$
Then, by a simple, but very long computation, we prove that
$((\Lambda,E),\omega)$ satisfies (\ref{eq-def-1}) and
(\ref{eq-def-2}). Thus, $((\Lambda,E),\omega)$ endows $M$ with a
twisted Jacobi structure.

\section{The characteristic foliation of a twisted Jacobi manifold}

It is well known \cite{dlm} that any Jacobi manifold is decomposed
into leaves equipped with transitive Jacobi structures that are,
according to the parity of the dimension of the leaves, contact or
locally conformal symplectic structures. In this section, we will
prove a similar result for twisted Jacobi manifolds.

\vspace{3mm}

Let $(M,(\Lambda,E),\omega)$ be a twisted Jacobi manifold and
consider its associated Lie algebroid over $M$, $(T^*M\times \R,
\{\cdot,\cdot \}^{\omega},\pi \circ (\Lambda,E)^\#)$. The image
$\mathrm{Im}(\pi \circ (\Lambda,E)^\#)$ of the anchor map defines
a completely integrable distribution on $M$, called the
\emph{characteristic distribution} of $(M,(\Lambda,E),\omega)$,
that determines a foliation of $M$ into leaves, which are called
the \emph{characteristic leaves of} $((\Lambda,E),\omega)$,
\cite{acw}. If, at every point of $M$, the dimension of the
characteristic leaf of $((\Lambda,E),\omega)$ through this point
is equal to the dimension of $M$, the twisted Jacobi manifold
$(M,(\Lambda,E),\omega)$ is said to be \emph{transitive}.
According to the parity of the dimension of $M$, there are two
kinds of transitive twisted Jacobi manifolds.

\begin{proposition} \label{p3.1}
Let $(M,(\Lambda,E),\omega)$ be a transitive twisted Jacobi
manifold.
\begin{enumerate}
\item[1)] If $M$ is of even dimension, then $((\Lambda,E),\omega)$
comes from a twisted locally conformal symplectic structure.
\item[2)] If $M$ is of odd dimension, then $((\Lambda,E),\omega)$
comes from a twisted contact structure.
\end{enumerate}
\end{proposition}
\begin{proof} 1) Let $\dim M=2n$. Since $((\Lambda,E),\omega)$
is transitive, $\mathrm{rank}\Lambda^{\#}=2n$, everywhere on $M$,
and $E$ is a section of $\mathrm{Im}\Lambda^{\#}$, i.e. there
exists a $1$-form $\vartheta$ on $M$ such that
$E=\Lambda^{\#}(\vartheta)$. Let $\Theta$ be the $2$-form on $M$
obtained by the inversion of $\Lambda$, i.e., for any $\alpha \in
\Gamma(T^*M)$, $i(\Lambda^{\#}(\alpha))\Theta = -\alpha$. A simple
computation shows that equations (\ref{eq-def-1}) and
(\ref{eq-def-2}) give, respectively, $d(\Theta - \omega) +
\vartheta \wedge (\Theta - \omega) =0$ and $d\vartheta =0$.
Whence, we conclude that $((\Lambda,E),\omega)$ is provided by the
twisted locally conformal symplectic structure
$(\vartheta,\Theta,\omega)$ on $M$.

\vspace{2mm}

\noindent 2) Let $\dim M=2n+1$. Since $((\Lambda,E),\omega)$ is
transitive, $\mathrm{rank}\Lambda^{\#}=2n$, everywhere on $M$, and
$E$ is not a section of $\mathrm{Im}\Lambda^{\#}$. Let $\vartheta$
be the $1$-form on $M$ defined by $i(E)\vartheta =1$ and
$\Lambda^{\#}(\vartheta)=0$ and let $\Theta$ be the $2$-form on
$M$ obtained by the inversion of $\Lambda$, i.e., for any $\alpha
\in \Gamma(T^*M)$, $i(\Lambda^{\#}(\alpha))\Theta = -(\alpha-
\langle \alpha,E\rangle \,\vartheta)$ and $i(E)\Theta =0$.
Clearly, $\vartheta\wedge \Theta^n \neq 0$, everywhere on $M$, and
$\Lambda = \Lambda^{\#}(\Theta)$. So, we have
$[\Lambda,\Lambda]=2\Lambda^{\#}(d\Theta)-2E\wedge
\Lambda^{\#}(d\vartheta)$ and, by a simple argumentation, we prove
that (\ref{eq-def-1}) and (\ref{eq-def-2}) give $\Theta =
d\vartheta + \omega$. Thus, $((\Lambda,E),\omega)$ comes from the
twisted contact structure $(\vartheta,\omega)$ on $M$.
\end{proof}

\begin{theorem} \label{t3.2}
Let $(M,(\Lambda,E),\omega)$ is a twisted Jacobi manifold. Then,
the bracket (\ref{croch-fun}) induces a transitive twisted Jacobi
structure on each  characteristic leaf of $M$.
\end{theorem}
\begin{proof}
Let $S$ be a characteristic leaf of $(M,(\Lambda,E),\omega)$
through a point $p$, with $\dim S =k$, and
$(x_1,\ldots,x_k,y_1,\ldots,y_{n-k})$, $n=\dim M$, a system of
adapted local coordinates of $M$. Given two functions
$\tilde{f},\tilde{g}\in C^{\infty}(S,\R)$, we can extend them
locally to functions $f,g \in \C$, i.e. $f(x,0)=\tilde{f}(x)$ and
$g(x,0)=\tilde{g}(x)$. On $C^{\infty}(S,\R)$ we define the bracket
$\{\,,\,\}_S$ by setting
\begin{equation}\label{br-S}
\{\tilde{f},\tilde{g}\}_S(x)=\{f,g\}(x,0), \hspace{7mm}
\mathrm{for}\;\;\mathrm{all}\;\; \tilde{f},\tilde{g} \in
C^{\infty}(S,\R).
\end{equation}
We have,
\begin{eqnarray*}
\{\tilde{f},\tilde{g}\}_S(x)=\{f,g\}(x,0) & = &
(\Lambda^{\#}(df)+fE)\vert_{(x,0)}g - \langle df,
E\rangle\vert_{(x,0)}g \\
& = & -(\Lambda^{\#}(dg)+gE)\vert_{(x,0)}f + \langle dg,
E\rangle\vert_{(x,0)}f  \\
\end{eqnarray*}
and we realize that the bracket (\ref{br-S}) only depends on
$\tilde{f}$ and $\tilde{g}$ because it is computed along the
integral curves of the vector fields $\Lambda^{\#}(df)+fE$,
$\Lambda^{\#}(dg)+gE$ and $E$ through $(x,0)$, which lie on $S$.
Clearly, (\ref{br-S}) yields a transitive twisted Jacobi structure
on $S$.
\end{proof}

\vspace{2mm}

From Proposition \ref{p3.1} and Theorem \ref{t3.2}, we conclude
that a twisted Jacobi manifold is foliated into leaves that are
endowed, according to the parity of the dimension, with a twisted
locally conformal symplectic structure or a twisted contact
structure.

\vspace{5mm}

\noindent{\bf Acknowledgments.} The work of Joana M. Nunes da
Costa has been partially supported by POCI/MAT/58452.


\begin{thebibliography}{99}

\bibitem{acw}
\textsc{A. Cannas da Silva and A. Weinstein}, \emph{Geometric
Models for Noncommutative Algebras}, University of California,
Berkeley Mathematics Lecture Notes 10 - AMS, Providence (1999).

\bibitem{dlm}
\textsc{P.~Dazord, A.~Lichnerowicz, C.-M.~Marle},  Structure
locale des vari\'et\'es de Jacobi, \emph{J. Math. Pures Appl.}
\textbf{70} (1991), 101.

\bibitem{im}
\textsc{D.~Iglesias and J.C.~Marrero}, Generalized Lie
bialgebroids and Jacobi structures, \emph{J. Geom. Phys.}
\textbf{40} (2001), 176.

\bibitem{krb}
\textsc{Y.~Kerbrat and Z.~Souici-Benhammadi},  Vari\'et\'es de
Jacobi et groupo\"{\i}des de contact, \emph{ C. R. Acad. Sci.
Paris, S\'erie I} \textbf{317} (1993), 81.


\bibitem{kl}
\textsc{C.~Klim\v{c}ic and T.~Strobl},  WZW-Poisson manifolds,
\emph{J. Geom. Phys.} \textbf{43} (2002), 341.


\bibitem{lch}
\textsc{A.~Lichnerowicz},  Les vari\'et\'es de Jacobi et leurs
alg\`ebres de Lie associ\'ees, \emph{J. Math. pures et appl.}
\textbf{57} (1978), 453.


\bibitem{jf}
\textsc{J.M.~Nunes da Costa and F.~Petalidou},  Twisted Jacobi
manifolds, \emph{Proceedings of XIV Fall Workshop on Geom. and
Phys. (Bilbao, 2005) Public. RSME}, vol. 8 (2006), 279.

\bibitem{jf-tj}
\textsc{J.M.~Nunes da Costa and F.~Petalidou}, Twisted Jacobi
manifolds, twisted Dirac-Jacobi structures and quasi-Jacobi
bialgebroids, \emph{J. Phys. A: Math. Gen.} \textbf{39} (2006),
10449.

\bibitem{p}
\textsc{J.S.-Park}, Topological open $p$-branes, in
\emph{Symplectic geometry and Mirror Symmetry} (Seoul, 2000), K.
Fukaya, Y.-G. Oh, K. Ono and G.Tian eds., World Sci. Publishing,
River Edge, NJ (2001), 311.

\bibitem{sw}
\textsc{P.~\v{S}evera and A.~Weinstein},  Poisson geometry with a
3-form background, \emph{Noncommutative geometry and string
theory, Prog. Theor. Phys.}, Suppl. \textbf{144} (2001), 145.

\end{thebibliography}
\end{document}